%

\long\def\comm(#1)#2!!!{}

{\gdef\editer
                                {\tolerance 3000
                                \gdef\\{{\tt\char'134}}\gdef\#{{\tt\char'043}}          
        
                                \gdef\{{{\tt\char'173}}\gdef\}{{\tt\char'175}}
                                \long\gdef\comm(##1)##2!!!
                                                        
        {{\leftskip20mm\parindent=0mm\medskip\hskip-20mm
                                                                \setbox1=\hbox{\\\tt##1}
                                                                {\ifdim\wd1<17mm\hbox to 
19mm{\box1\hfill}%
                                                                \else\hbox to \wd1{\unhbox1} : 
\fi}##2.\par}}
                                \centerline{\bfXII Les commandes de DehTeX}\bigskip
                                \datef\bigskip\bigskip
                                Num\'eros des familles de fontes: 
        0=cmr$\fl$rm; 1=cmmi$\fl$mi; 2=cmsy$\fl$sy; 3=cmex;  4=cmti$\fl$it;
                                5=cmsl$\fl$sl; 6=cmbx$\fl$bf; 7=cmtt$\fl$tt; 8=cmb$\fl$be; 
                                9=cmmib$\fl$bm; A=msam$\fl$sa; B=msbm$\fl$sb; 
C=eufm$\fl$go;
                                D=eurm$\fl$eu; E non utilis\'e; F non utilis\'e;
                                \def\dump{}}}


\def\<<{{\cyr\char'074~}}
                                \comm(<<)Ecrit \<<(avec blanc inseccable apr\`es)!!!

\def\>>{~{\cyr\char'076}}
                                \comm(>>)Ecrit \>> (avec blanc inseccable avant)!!!

\def\\{\hbox{$\backslash$}}
                                \comm($\\$)Ecrit $\\$; ne n\'ecessite pas de blanc apr\`es!!! 

\let\(=\langle
                                \comm(()Ecrit $\($; mode math\'ematique!!!

\let\)=\rangle
                                \comm({)})Ecrit $\)$; mode math\'ematique!!!

\def\0{{\mi 0}}\def\1{{\mi 1}}\def\2{{\mi 2}}\def\3{{\mi 3}}
\def\4{{\mi 4}}\def\5{{\mi 5}}\def\6{{\mi 6}}\def\7{{\mi 7}}
\def\8{{\mi 8}}\def\9{{\mi 9}}
                                \comm({\it chiffre})Ecrit le chiffre correspondant en 
\<<oldstyle\>>: \0,
                                                                \1, \2, \3, \4, \5, \6, \7, \8, \9!!!

\font\V=cmr5 \font\VI=cmr6 \font\VII=cmr7 \font\VIII=cmr8
\font\IX=cmr9 \font\XII=cmr10 scaled\magstep1
\font\XIV=cmr10 scaled\magstep2
\font\XVIII=cmr10 scaled\magstep3
\font\XXIV=cmr10 scaled\magstep4
                                \comm({{\it chiffre romain}})Appelle la taille de caract\`ere
                                correspondante dans la fonte cmr; disponibles: {\tt\\V, \\VI, \\VII,
                                \\VIII, \\IX, \\XII, \\XIV, \\XVIII, \\XXIV}!!!

\mathchardef\a="710B
                                \comm(a)Ecrit $\a$; mode math\'ematique!!!

\def\adots{\mathinner{\mkern2mu\raise1pt\hbox{.}
                                                                \mkern3mu\raise4pt\hbox{.}
                                                                \mkern1mu\raise7pt\hbox{.}}}
                                \comm(adots)Ecrit $\adots$; mode math\'ematique!!!

\let\al=\aleph
                                \comm(al)Ecrit $\al$; mode math\'ematique!!!

\def\and{\hbox{ and }}
                                \comm(and)Ecrit $\and$ en romain avec blanc avant et apr\`es, 
quel
                                                                que soit le mode (le m\^eme en 
pench\'e s'appelle {\tt\\sland})!!!

\def\ape{\mathchar"3A44}
                                \comm(ape)Ecrit $\ape$; mode math\'ematique!!!
                 
\def\app{\mathord{\in}}
                                \comm(app)Ecrit $\app$ sans blanc avant ni apr\`es; mode 
math\'ematique!!!

\let\aps=\triangleright
                                \comm(aps)Ecrit $\aps$; mode math\'ematique!!!
                
\let\av=\triangleleft
                                \comm(av)Ecrit $\av$; mode math\'ematique; autre nom: {\tt 
avs}!!!
                
\def\ave{\mathchar"3A45}
                                \comm(ave)Ecrit $\ave$; mode math\'ematique!!!

                                \comm(avs)Voir {\tt ave}!!!

\mathchardef\b="710C       
                                \comm(b)Ecrit $\b$; mode math\'ematique!!!

\def\bbar#1{\overline{\mathstrut#1}}
                                \comm(bbar)La commande {\tt\\bbar\#} place une barre sur \#, \`a
                                                                une hauteur uniforme et assez 
grande; mode math\'ematique!!!

                                                                \font\cmbX=cmb10 at 10pt 
\font\cmbVII=cmb10 at 7pt 
                                                                \font\cmbV=cmb10 at 5pt 
\textfont8=\cmbX
                                                        
        \scriptfont8=\cmbVII\scriptscriptfont8=\cmbV
                                                                \def\be{\fam8\cmbX}
                                                                \let\bfe=\be
                                \comm(be)Passe \`a la famille \<<gras \'etroit\>> (fontes cmb); 
{\be
                                                                Ceci est du be}; autre nom: 
{\tt\\bfe}!!!

                                                        
        \font\cmbxX=cmbx10\font\cmbxVII=cmbx7 
                                                                \font\cmbxV=cmbx5 
\textfont6=\cmbxX
                                                        
        \scriptfont6=\cmbxVII\scriptscriptfont6=\cmbxV
                                                                \def\bf{\fam6\cmbxX}

\font\bfXII=cmbx10 scaled\magstep1  
              \font\bfXIV=cmbx10 scaled\magstep2  
              \font\bfXVIII=cmbx10 scaled\magstep3  
              \font\bfXXIV=cmbx10 scaled\magstep4
                                \comm(bf{{\it chiffre romain}})Passe \`a la fonte \<<gras\>>\ 
(cmbx)
                                                                de la taille correspondante;
                                                                disponibles: {\tt\\bfXII, \\bfXIV, 
\\bfXVIII, \\bfXXIV}!!!

\def\bg{\medbreak\medskip}
        \comm(bg)Provoque un grand saut et favorise la coupure; la moiti\'e du
                                saut est supprim\'ee apr\`es un display!!!

\def\bgni{\medbreak\medskip\noindent}
                                \comm(bgni)Provoque un grand saut, favorise la coupure et
                                supprime l'indentation apr\`es!!!

\let\bi=\cmbxtiX
                                \comm(bi)Passe \`a la fonte \<<gras italique\>> (fontes cmbxti); 
{\bi
                                Ceci est du bi}!!!

                                \comm(bl)Ecrit $\bullet$!!!

                                                        
        \font\cmmibX=cmmib10\font\cmmibVII=cmmib7
                \font\cmmibV=cmmib5 \textfont9=\cmmibX
                \scriptfont9=\cmmibVII\scriptscriptfont9=\cmmibV
                \def\bm{\fam9\cmmibX}
                                \comm(bm)Appelle la famille \<<gras math\'ematique\>> (fontes
                                                                cmmib); {\bm Ceci est du bm}!!!

\def\bn{\bigskip\nobreak}
                                \comm(bn)Provoque un saut et d\'efavorise la coupure!!!

                                \comm(bnni)Provoque un saut, d\'efavorise la coupure et 
supprime
                                                                l'indentation!!!

\long\def\boxit#1#2{\setbox1=\hbox{\kern#1{#2}\kern#1}%
                                                                \dimen1=\ht1 \advance\dimen1 by 
#1
                                                                \dimen2=\dp1 \advance\dimen2 by 
#1
                                                                \setbox1=\hbox{\vrule 
height\dimen1 depth\dimen2\box1\vrule}%
                                                        
        \setbox1=\vbox{\hrule\box1\hrule}%
                                                                \advance\dimen1 by .4truept 
\ht1=\dimen1
                                                                \advance\dimen2 by .4truept 
\dp1=\dimen2 \box1\relax}
                                \comm(boxit)Encadre son second argument; syntaxe
                                                                {\tt\\boxit$\{$2pt$\}\{$xx$\}$} 
donne \boxit{2pt}{{\tt xx}}!!!

\let\bs=\cmbxslX
                                \comm(bs)Passe \`a la fonte \<<gras pench\'e\>> (fontes cmbxsl); 
{\bs
                                ceci est du bs}!!!

\def\card{{\rm card}}
                                \comm(card)Ecrit \card!!!               

\def\carre{\hbox{\sa\char'004}}
                                \comm(carre)Ecrit \carre!!!

\catcode`\@=11\def\cases#1{\left\{\,\vcenter{\normalbaselines\m@th
                                                        
        \ialign{$##\hfil$&\quad{\rm##}\hfil\crcr#1\crcr}}\right.}
                                                                \catcode`\@=12

\def\cf{\hbox{\it c.f. }}
                                \comm(cf)Ecrit \cf\ avec un espace apr\`es!!!

\font\cmdunhX=cmdunh10
                                \comm(cmdunhX)Passe \`a la fonte cmdunh10; {\cmdunhX Ceci 
est du
                                                                cmdunh10}!!!

\def\comp{{\scriptscriptstyle\circ}}
                                \comm(comp)Ecrit $\comp$; mode math\'ematique!!!

                                \comm(Cor)La commande {\tt Cor\#} provoque un saut, puis 
\'ecrit
                                                                sans indentation {\bf Corollary\#.-
}!!!

                                \comm(cyr)Fait passer dans la fonte \<<wncyr10\>>; {\cyr Ceci 
est du
                                                                cyr}!!!

\mathchardef\d="710E
                                \comm(d)Ecrit $\d$; mode math\'ematique!!!

\mathchardef\D="7101
                                \comm(D)Ecrit $\D$; mode math\'ematique!!!

\def\date{\par\line{\hfil\ifcase\month\or January\or
                                                                February\or March\or April\or 
May\or June\or July\or August\or
                                                                September\or October\or 
November\or December\fi\ 
                                                                {\oldstyle\the\day},\ 
{\oldstyle\the\year}}}
                                \comm(date)Met la date en anglais \`a droite de la ligne!!!

\def\datef{\par\line{\hfil le\ {\oldstyle\the\day}\
                 \ifcase\month\or
                 janvier\or f\'evrier\or mars\or avril\or mai\or juin\or
                 juillet\or ao\^ut\or septembre\or octobre\or novembre\or
                 d\'ecembre\fi\ {\oldstyle\the\year}}}
                                \comm(datef)Met la date en francais \`a droite de la ligne!!!

                                \comm(Def)Provoque un saut, puis \'ecrit {\bf Definition.-}!!!

\let\der=\partial
                                \comm(der)Ecrit $\der$; mode math\'ematique!!!

                                \comm(diagram)Pour coder un diagramme; m\^eme syntaxe que
                                                                {\tt\\matrix}, simplement les espaces 
sont red\'efinis!!!

\def\Dom{\mathord{\rm Dom}}
                                \comm(Dom)Ecrit $\Dom$; mode math\'ematique!!!

\mathchardef\e="7122
                                \comm(e)Ecrit $\e$; mode math\'ematique!!!

\def\eg{\hbox{\it e.g. }}
                                \comm(eg)Ecrit \eg avec un espace apr\`es!!!

                                \comm(entete)Place l'entete de Caen en haut de page; attention! il 
faut
                                                                que le sceau soit dans le fichier 
pictures!!!

\let\equ=\Longleftrightarrow
                                \comm(equ)Ecrit $\equ$; mode math\'ematique!!!

                                \comm(esp)Force un blanc, m\^eme dans les fontes o\`u l'espace 
est nul,
                                                                comme cmmi qui sert dans 
oldstyle!!!

\def\et{\hbox{ et }}
                                \comm(et)Ecrit $\et$ en tout mode (avec un espace avant et 
apr\`es)!!!

\def\etc{\hbox{\it etc\dots\ }}
                                \comm(etc)Ecrit \etc avec un espace apr\`es!!!

                                                        
        \font\euX=eurm10\font\euVII=eurm7\font\euV=eurm5
                                                                \textfont13=\euX
                                                        
        \scriptfont13=\euVII\scriptscriptfont13=\euV
                                                                \def\eu{\fam13\euX}
                                \comm(eu)Passe \`a la famille \<<euler\>> (fontes eurm); {\eu
                                                                Ceci est du eu}!!!

\mathchardef\eup="0D3A
                                \comm(eup)Ecrit $\eup$ (point en fonte \<<euler\>>); mode
                                                                math\'ematique!!!

\mathchardef\euv="0D3B
                                \comm(euv)Ecrit $\euv$ (virgule en fonte \<<euler\>>); mode
                                                                math\'ematique!!!

\let\ev=\emptyset
                                \comm(ev)Ecrit $\ev$; mode math\'ematique!!!

\let\ex=\exists
                                \comm(ex)Ecrit $\ex$; mode math\'ematique!!!

                                \comm(Ex)Provoque un grand saut, puis \'ecrit {\bf Example.-
}!!!

\mathchardef\f="7127
                                \comm(f)Ecrit $\f$; mode math\'ematique!!!

\mathchardef\F="7108
                                \comm(F)Ecrit $\F$; mode math\'ematique!!!
                 
\let\fl=\longrightarrow
                                \comm(fl)Ecrit $\fl$; mode math\'ematique!!!

                                \comm(Fin)Ecrit $\carre$ pr\'ec\'ed\'e d'un blanc inseccable!!!

\let\flc=\rightarrow
                                \comm(flc)Ecrit $\flc$ (fl\`eche courte); mode math\'ematique!!!

{\catcode`@=11
                                                                \catcode`\;=\active%
                                                                \catcode`\:=\active%
                                                                \catcode`\!=\active%
                                                                \catcode`\?=\active%
                                \gdef\francais{
                                                                \catcode`\;=\active%
                                                        
        \def;{\relax\ifhmode\ifdim\lastskip>\z@\unskip\fi
                                                                \kern.2em\fi\string;}
                                                                \catcode`\:=\active%
                                                        
        \def:{\relax\ifhmode\ifdim\lastskip>\z@\unskip\fi
                                                                \kern.2em\fi\string:}
                                                                \catcode`\!=\active%
                                                        
        \def!{\relax\ifhmode\ifdim\lastskip>\z@\unskip\fi
                                                                \kern.2em\fi\string!}
                                                                \catcode`\?=\active%
                                                        
        \def?{\relax\ifhmode\ifdim\lastskip>\z@\unskip\fi
                                                                \kern.2em\fi\string?}
                                                                \frenchspacing\catcode`\@=12}}
                                \comm(francais)Adapte \`a la typographie francaise;
                                                                redefinit les caracteres ; : ? !  aux 
normes francaises, c'est \`a
                                                                dire avec un blanc inseccable devant; 
cela ne change rien de taper
                                                                ou de ne pas taper les blancs avant 
les caracteres ; : ! ? (par contre
                                                                cela change d'en taper apres ou 
non)!!!

\mathchardef\g="710D
                                \comm(g)Ecrit $\g$; mode math\'ematique!!!

\mathchardef\G="7100
                                \comm(G)Ecrit $\G$; mode math\'ematique!!!
                  
                                                        
        \font\eufmX=eufm10\font\eufmVII=eufm7\font\eufmV=eufm5
                                                        
        \textfont12=\eufmX\scriptfont12=\eufmVII
                                                                \scriptscriptfont12=\eufmV
                                                                \def\go{\fam12\eufmX}
                                \comm(go)Passe dans la famille \<<gothique\>> (fontes eufm); 
                                                                {\go Ceci est du go}!!!

\mathchardef\gpref="3A40
                                \comm(gpref)Ecrit $\gpref$ (le g signifie \<<grand\>>) ; mode
                                                                math\'ematique!!!
                                                        
\let\gprefe=\sqsubseteq
                                \comm(gprefe)Ecrit $\gprefe$ (le g signifie \<<grand\>>) ; mode 
                                                                math\'ematique!!!

\mathchardef\gsuff="3A41
                                \comm(gsuff)Ecrit $\gsuff$ (le g signifie \<<grand\>>) ; mode
                                                                math\'ematique!!!
                                                        
\let\gsuffe=\sqsupseteq
                                \comm(gsuffe)Ecrit $\gsuffe$ (le g signifie \<<grand\>>) ; mode 
                                                                math\'ematique!!!

\mathchardef\h="7112
                                \comm(h)Ecrit $\h$; mode math\'ematique!!!

\mathchardef\H="7102
                                \comm(H)Ecrit $\H$; mode math\'ematique!!!
                  
\font\hel=cmb10 at 10 pt
               \font\helVII=cmb10 at 7 pt
               \font\helXII=cmb10 at 12 pt
               \font\helXIV=cmb10 at 14 pt
                                \comm(hel)Passe dans la fonte \<<cmb10\>>; 
                                                                {\hel Ceci est du hel};
                                                                suivie \'eventuellemnt d'une taille en 
romain majuscule;
                                                                disponibles: {\tt \\helVII, \\helXII, 
\\helXIV}!!!

\font\helbf=cmb10 at 10 pt
               \font\helbfXII=cmb10 at 12 pt
               \font\helbfXIV=cmb10 at 14 pt
               \font\helbfXVIII=cmb10 at 18 pt
               \font\helbfXXIV=cmb10 at 24 pt
               \font\helbfXXXVI=cmb10 at 36 pt
                                \comm(helbf)Passe \`a la fonte \<<cmb10 gras\>>;
                                                         {\helbf ceci est du helbf};
                                                                suivie \'eventuellemnt d'une taille en 
chiffre romain majuscule; 
                                                                disponibles:
                                                                {\tt\\helbfXII, \\helbfXIV, 
\\helbfXVIII,\\helbfXXIV,
                                                                \\helbfXXXVI}!!!

\font\helbi=cmb10 at 10 pt
               \font\helbiXII=cmb10 at 12 pt
               \font\helbiXIV=cmb10 at 14 pt
               \font\helbiXVIII=cmb10 at 18 pt
               \font\helbiXXIV=cmb10 at 24 pt
               \font\helbiXXXVI=cmb10 at 36 pt
                                \comm(helbi)Passe \`a la fonte \<<cmb10 gras italique\>>; 
                                                                {\helbi Ceci est du helbi}; 
                                                                suivie \'eventuellemnt d'une taille en 
chiffres romains majuscules;
                                                                disponibles:
                                                                {\tt\\helbiXII, \\helbiXIV, 
\\helbiXVIII, \\helbiXXIV, 
                                                                \\helbiXXXVI}!!!

\def\hfl#1#2{\smash{\mathop{\hbox to 12truemm{\rightarrowfill}}
                                                        
        \limits^{\scriptstyle#1}_{\scriptstyle#2}}}
                                \comm(hfl{\#1}{\#2})Trace une fl\`eche horizontale de 12 mm 
avec
                                                                {\tt\#1} dessus et {\tt\#2} dessous!!!

                                \comm(hhat)Ecrit un grand chapeau!!!

\def\ie{\hbox{\it i.e. }}
                                \comm(ie)Ecrit \ie avec un blanc apr\`es!!!

\def\iff{if and only if }
                                \comm(iff)Ecrit \iff avec un blanc apr\`es!!!

\def\Im{{\be Im}}
                                \comm(Im)Ecrit \Im!!!

\let\imp=\Longrightarrow
                                \comm(imp)Ecrit $\imp$; mode math\'ematique; autre nom:
                                                                {\tt\\impl}!!!

\let\impl=\Longrightarrow
                                \comm(impl)Voir {\tt\\imp}!!!

\let\inc=\subset
                                \comm(inc)Ecrit $\inc$; mode math\'ematique!!!

\let\ince=\subseteq
                                \comm(ince)Ecrit $\ince$; mode math\'ematique!!!

\def\inserer(#1)(#2)(#3)(#4)(#5){
                                                                \dimen1=#2\divide\dimen1 
by1000\multiply\dimen1 by#4
                                                                \dimen2=#3\divide\dimen2 
by1000\multiply\dimen2 by#4
                                                        
        \midinsert\vglue\dimen2$$\vbox{\hsize=\dimen1
                                                                \ni\special{picture #1 scaled #4}}$$
                                                                \centerline{#5}\endinsert}
                                \comm(inserer)Ins\`ere au mieux un dessin
                                                                la syntaxe est {\tt\\inserer(nom de la 
figure)(largeur)(hauteur)
                                                                (echelle)(commentaire)}
                                                                (les parametres largeur et hauteur
                                                                sont \`a lire dans le fichier Pictures; 
il faut taper l'unite)!!!

\let\inter=\cap
                                \comm(inter)Ecrit $\inter$; mode math\'ematique!!!

\let\Inter=\bigcap
                                \comm(Inter)Ecrit $\Inter$; mode math\'ematique!!!

                                                        
        \font\cmtiX=cmti10\font\cmtiVII=cmti7 
                                                                \font\cmtiV=cmti10 at 5pt 
                                                        
        \textfont4=\cmtiX\scriptfont4=\cmtiVII\scriptscriptfont4=\cmtiV
                                                                \def\it{\fam4\cmtiX}

\mathchardef\k="7114
                                \comm(k)Ecrit $\k$; mode math\'ematique!!!

\mathchardef\l="7115
                                \comm(l)Ecrit $\l$; mode math\'ematique!!!

\mathchardef\L="7103
                                \comm(L)Ecrit $\L$; mode math\'ematique!!!
                  
\def\Lem#1{\bgni{\bf Lemma#1.-}}
                                \comm(Lem)La commande {\tt\\Lem\#} provoque un saut, puis 
\'ecrit
                                                                sans indentation {\bf Lemma\#.-}!!!

                                \comm(lettre)Ecrit le haut de page pour une lettre en anglais; il faut
                                                                le sceau dans le fichier pictures!!!

                                \comm(lettref)Ecrit le haut de page pour une lettre en francais; il 
faut
                                                                le sceau dans le fichier pictures!!!

                                \comm(lettrem)Ecrit le haut de page pour une lettre en francais!!!

\def\lpol{\leavevmode\setbox0=\hbox{l}\hbox to \wd0{\hss\char'40l}}
                                \comm(lpol)Ecrit \lpol\ (l polonais)!!!

\def\Lpol{\leavevmode\setbox0=\hbox{L}\hbox to \wd0{\hss\char'40L}}
                                \comm(Lpol)Ecrit \Lpol\ (L polonais)!!!

\mathchardef\m="7116
                                \comm(m)Ecrit $\m$; mode math\'ematique!!!

                                \comm(mg)Provoque un saut moyen et favorise la coupure!!!

                                \comm(mgni)Provoque un saut moyen, favorise la coupure et 
supprime
                                                                l'indentation!!!

                                                                \def\mi{\fam1\teni}
                                \comm(mi)Passe \`a la famille \<<math. italique\>> (fontes cmmi); 
{\mi
                                                                ceci est du mi}; autres noms: 
{\tt\\mit} (mode math\'ematique),
                                                                {\tt\\oldstyle} (tout mode)!!!

                                \comm(mn)Provoque un saut moyen et d\'efavorise la coupure!!!

                                \comm(mnni)Provoque un saut moyen, d\'efavorise la coupure et
                                                                supprime l'indentation!!!

\def\move#1#2#3
                                                                {\vbox to0pt{\kern-
#3mm\smash{\hbox{\kern#2mm{#1}}}\vss}
                                                                \nointerlineskip}
                                \comm(move\#1\#2\#3)Ecrit \#1 \`a une abscisse de \#2 mm et \`a
                                                                une ordonn\'ee de \#3 mm par 
rapport au point courant (qui n'est
                                                                pas modifi\'e); mode vertical. 
Utilisable pour figures. Fabriquer
                                                                une \\{\tt vbox} suivant la syntaxe 
suivante \par
                                                                \hskip3cm{\tt \\vbox to ... 
mm\{\\hsize=... mm\\vfill\par
                                                         \hskip3cm\\special \{picture ... scaled ... \} 
\par 
                                                         \hskip3cm\\move \{\ ... \}\{\ ... \}\{\ ... \} 
\par 
                                                                \hskip3cm\\hfill\} }\par
                                                                Le {\tt\\vfill} est pour tasser la boite 
vers le bas, le {\tt\\hfill}
                                                                est pour tasser vers la gauche et 
assurer la largeur!!!

\mathchardef\n="7117
                                \comm(n)Ecrit $\n$; mode math\'ematique!!!

\def\Nat{\hbox{\rm I\kern-.9truemm\hbox{N}}}
                                \comm(Nat)Ecrit \Nat!!!

\let\ni=\noindent
                                \comm(ni)Supprime l'indentation!!!

\mathchardef\NN="0B4E
                                \comm(NN)Ecrit $\NN$: mode math\'ematique!!!

\def\non{\mathord{\be non}}
                                \comm(non)Ecrit $\non$; mode math\'ematique!!!

                                \comm(notapp)Ecrit $\notin$ sans blanc avant ni apr\`es; mode
                                                                        math\'ematique!!!

\mathchardef\o="7121
                                \comm(o)Ecrit $\o$; mode math\'ematique!!!

\mathchardef\O="710A
                                \comm(O)Ecrit $\O$; mode math\'ematique!!!
                  
\let\oo=\infty
                                \comm(oo)Ecrit $\oo$; mode math\'ematique!!!

\def\ou{\hbox{ ou }}
                                \comm(ou)Ecrit $\ou$ en tout mode (avec un blanc avant et 
apr\`es)!!!

\mathchardef\p="7119
                                \comm(p)Ecrit $\p$; mode math\'ematique!!!

\mathchardef\P="7105
                                \comm(P)Ecrit $\P$; mode math\'ematique!!!

\def\pmb#1{\setbox0=\hbox{#1}\kern-.15pt\copy0\kern-\wd0
                 \kern-.3pt\copy0\kern-\wd0\kern-.15pt\raise.1pt\box0}
                                \comm(pmb)Met en gras son param\`etre: {\tt\\pmb\#} \'ecrit \# en
                                                                faux gras (ce qui donne 
\pmb{\#})!!!

\def\pp{\hbox{$\ldots$}}
                                \comm(pp)Ecrit \pp!!!

\def\Ppar{\mathhexbox27B}
                                \comm(Ppar)Ecrit \Ppar!!!

\def\ppp{\hbox{$, \ldots, $}}
                                \comm(ppp)Ecrit \ppp!!!

\let\pr=\vdash
                                \comm(pr)Ecrit $\pr$; mode math\'ematique!!!

\def\Pr{\bgni{\it Proof. }}
                                \comm(Pr)Fait un saut et \'ecrit \<<{\it Proof.}\>>!!!

\def\pref{\mathrel{\scriptstyle\gpref}}
                                \comm(pref)Ecrit la relation $\pref$; mode math\'ematique!!!
                                                        
\def\prefe{\mathrel{\scriptstyle\gprefe}}
                                \comm(prefe)Ecrit la relation $\prefe$; mode math\'ematique!!!
                                                        
\def\Prop#1{\bgni{\bf Proposition#1.-}}
                                \comm(Prop)La commande {\tt\\Prop\#} provoque un saut, puis 
\'ecrit
                                                                sans indentation {\bf Proposition\#.-
}!!!

\mathchardef\Pw="0C50
                                \comm(Pw)Ecrit $\Pw$; mode math\'ematique!!!

\mathchardef\q="711F
                                \comm(q)Ecrit $\q$; mode math\'ematique!!!

\let\qq=\forall
                                \comm(qq)Ecrit $\qq$; mode math\'ematique!!!

\mathchardef\r="711A
                                \comm(r)Ecrit $\r$; mode math\'ematique!!!

\mathchardef\res="0A16
                                \comm(res)Ecrit $\res$; mode math\'ematique!!!
                                                        
\def\resp{\hbox{\it resp. }}
                                \comm(resp)Ecrit \resp (avec un blanc apr\`es)!!!

\def\Ree{\hbox{\rm I\kern-.9truemm\hbox{R}}}
                                \comm(Ree)Ecrit \Ree!!!

\mathchardef\RR="0B52
                                \comm(RR)Ecrit $\RR$; mode math\'ematique!!!

\def\rres{\hbox to 4.5pt{$\res$\kern-1pt\hbox{$\res$}\hss}}
                                \comm(rres)Ecrit $\rres$; mode math\'ematique!!!

\mathchardef\s="711B
                                \comm(s)Ecrit $\s$; mode math\'ematique!!!

\mathchardef\S="7106
                                \comm(S)Ecrit $\S$; mode math\'ematique!!!

                                                        
        \font\msamX=msam10\font\msamVII=msam7 
                                                                \font\msamV=msam5 
\textfont10=\msamX
                                                        
        \scriptfont10=\msamVII\scriptscriptfont10=\msamV
                                                                \def\sa{\fam10\msamX}
                                \comm(sa)Passe \`a la famille \<<symboles a\>> (fontes msam); 
{\sa
                                                                Ceci est du sa}!!!

                                                        
        \font\msbmX=msbm10\font\msbmVII=msbm7 
                                                                \font\msbmV=msbm5 
\textfont11=\msbmX
                                                        
        \scriptfont11=\msbmVII\scriptscriptfont11=\msbmV
                                                                \def\sb{\fam11\msbmX}
                                                                \let\db=\sb
                                \comm(sb)Passe \`a la famille \<<symboles b\>> (fontes msbm); 
{\sb
                                                                CECI EST DU SB}; autre nom: 
{\tt\\db}!!!

\let\sat=\models
                                \comm(sat)Ecrit $\sat$; mode math\'ematique!!!

                                \comm(sg)Provoque un saut petit et favorise la coupure!!!

\def\sgni{\smallbreak\noindent}
                                \comm(sgni)Provoque un saut petit, favorise la coupure et 
supprime
                                                                l'indentation!!!

                                                        
        \font\cmslX=cmsl10\font\cmslVII=cmsl10 at 7 pt 
                                                                \font\cmslV=cmsl10 at 5pt 
                                                        
        \textfont5=\cmslX\scriptfont5=\cmslVII\scriptscriptfont5=\cmslV
                                                                \def\sl{\fam5\cmslX}

\def\sland{\hbox{ \sl and }}
                                \comm(sland)Ecrit $\sland$ en tout mode!!!

\font\smcap=cmcsc10                                                             
                                \comm(smcap)Passe \`a la fonte \<<petites capitales\>>\ (cmcsc);
                                                                {\smcap Ceci est du smcap}!!!

                                \comm(sn)Provoque un saut petit et d\'efavorise la coupure!!!

                                \comm(snni)Provoque un saut petit, d\'efavorise la coupure et
                                                                supprime l'indentation!!!

\def\Spar{\mathhexbox278}
                                \comm(Spar)Ecrit \Spar!!!

\mathchardef\SS="0B53
                                \comm(SS)Ecrit $\SS$: mode math\'ematique!!!

\def\ssi{si et seulement si }
                                \comm(ssi)Ecrit \ssi (avec un blanc apr\`es)!!!

\def\suff{\mathrel{\scriptstyle\gsuff}}
                                \comm(suff)Ecrit la relation $\suff$; mode math\'ematique!!!
                                                        
\def\suffe{\mathrel{\scriptstyle\gsuffe}}
                                \comm(suffe)Ecrit la relation $\suffe$; mode math\'ematique!!!
                                                        
\def\Supp{\mathord{\be Supp}}
                                \comm(Supp)Ecrit $\Supp$; mode math\'ematique!!!

                                                                \def\sy{\fam2\tensy}
                                                                \let\ca=\sy
                                \comm(sy)Passe \`a la famille \<<symbole\>> (fontes cmsy); {\sy
                                                                CECI EST DU SY}; autres noms: 
{\tt\\cal} (mode math\'ematique),
                                                                {\tt\\ca} (tout mode)!!!

\mathchardef\t="711C
                                \comm(t)Ecrit $\t$; mode math\'ematique!!!

\def\Thm#1{\bgni{\bf Theorem#1.-}}
                                \comm(Thm)La commande {\tt\\Thm\#} provoque un saut, puis 
\'ecrit
                                                                sans indentation {\bf Theorem\#.-
}!!!

\font\tim=cmb10 at 10pt 
                
               \font\timXIV=cmb10 at 14pt 
               \font\timXVIII=cmb10 at 18pt 
               \font\timVII=cmb10 at 7pt 
                                \comm(tim)Passe dans la fonte \og cmb10\fg; 
                                                                {\tim Ceci est du tim};
                                                                suivie \'eventuellemnt d'une taille en 
chiffres romains majuscules;
                                                                disponibles: {\tt\\timVII, \\timlXII, 
\\timXIV, \\timXVIII}!!!

\font\timbf=cmb10 at 10pt 
               
              \font\timbfXIV=cmb10 at 14pt 
              \font\timbfXVIII=cmb10 at 18pt 
              \font\timbfVII=cmb10 at 7pt
                                \comm(timbf)Passe dans la fonte \og cmb10 gras\fg; 
                                                                {\timbf Ceci est du timbf};
                                                                suivie \'eventuellemnt d'une taille en 
chiffres romains majuscules;
                                                                disponibles: {\tt\\timbfVII, 
\\timbflXII, \\timbfXIV,
                                                                \\timbfXVIII}!!!

                                                        
        \font\cmttX=cmtt10\font\cmttVII=cmtt10 at 7 pt 
                                                                \font\cmttV=cmtt10 at 5pt 
                                                        
        \textfont7=\cmtiX\scriptfont7=\cmttVII\scriptscriptfont7=\cmttV
                                                                \def\tt{\fam7\cmttX}

                                \comm(ttil)Ecrit un grand tilde!!!

\mathchardef\u="711D
                                \comm(u)Ecrit $\u$; mode math\'ematique!!!

\mathchardef\U="7107
                                \comm(U)Ecrit $\U$; mode math\'ematique!!!

\let\un=\cup
                                \comm(un)Ecrit $\un$; mode math\'ematique!!!

\let\Un=\bigcup
                                \comm(Un)Ecrit $\Un$; mode math\'ematique!!!

\def\val{\mathop{\be val}}
                                \comm(val)Ecrit $\val$; mode math\'ematique!!!

\def\var{{\be var}}
                                \comm(var)Ecrit $\var$!!!

                                \comm(vfl{\#1}{\#2})Trace une fl\`eche verticale de 12 mm avec
                                                                {\tt\#1} \`a gauche et {\tt\#2} \`a 
droite!!!

\mathchardef\w="7120
                                \comm(w)Ecrit $\w$; mode math\'ematique!!!

\mathchardef\W="7109
                                \comm(W)Ecrit $\W$; mode math\'ematique!!!

                                \comm(wrt)Ecrit \<<with respect to\>> suivi d'un blanc 
inseccable!!!

\mathchardef\x="7118
                                \comm(x)Ecrit $\x$; mode math\'ematique!!!

\mathchardef\X="7104
                                \comm(X)Ecrit $\X$; mode math\'ematique!!!

\mathchardef\y="7111
                                \comm(y)Ecrit $\y$; mode math\'ematique!!!

\mathchardef\z="7110
                                \comm(z)Ecrit $\z$; mode math\'ematique!!!

\def\ZZ{{\db Z}}
                                \comm(ZZ)Ecrit $\ZZ$!!!

                                \comm( )\vfill\eject\centerline{\bfXII Anciennes commandes}!!!

\let\bfe=\be
                                \comm(bfe)Ancien nom de {\tt\\be}!!!

                                \comm(bfi)Ancien nom de {\tt\\bi}!!!
  
\let\ca=\sy
                                \comm(ca)Ancien nom de la famille {\tt\\sy}!!!

\font\cmmibX=cmmib10
                                \comm(cmmibX)Passe \`a la fonte cmmib10; {\cmmibX Ceci est 
du
                                                                cmmib}!!!

\let\db=\sb
                                \comm(db)Ancien nom de la famille {\tt\\sb}!!!

\def\fg{~{\cyr\char'076}}
                                \comm(fg)Ancien nom de {\tt\\>>}!!!

                                \comm(msymX)Passe \`a la fonte \<<msbm10\>>; obsol\`ete: 
nouveau
                                                                nom de famille: {\tt\\db}!!!

\def\og{{\cyr\char'074}~}
                                \comm(og)Ancien nom de {\tt\\<<}!!!

                                \comm(rd)Ancien nom de la famille {\tt\\mi}!!!

\def\picture #1 by #2 (#3){
 \vbox to #2{\hrule width #1 height 0pt depth 0pt
 \vfill\special{picture #3}}}
\def\dessin(#1)(#2)(#3)(#4){{
 \dimen0=#2 \dimen1=#3
 \divide\dimen0 by 1000 \multiply\dimen0 by #4
 \divide\dimen1 by 1000 \multiply\dimen1 by #4
 \picture \dimen0 by \dimen1 (#1 scaled #4)}}

\newcount\numlem
\newcount\numchap
\hsize=12cm
\vsize=20truecm
\voffset=1cm
\hoffset=5mm
\newtoks\gauche \gauche={Serge Burckel}
\newtoks\droite 
\def\makeheadline{\vbox to 0pt{\vskip -40pt\line{\vbox to 8.5pt{}\the
\headline}\vss\nointerlineskip}}
\headline={{\voffset 2cm\sevenbf\the\gauche\hfill\the\droite}}

\def\sk{\smallskip}

\overfullrule=0mm 
\parindent=0mm

 

\def\ca{{\bi a}}

\def\l{n}

\def\psi{{\bf cl}}
\def\x{X}
\def\fl{{\longrightarrow}}

\def\bar#1{{\overline{#1}}}

\def\H#1{{\hskip #1truemm}}
\def\BX#1{\boxit{8pt}{\vbox{#1}}}

\def\s#1 #2{{#1}_{(#2)}}
\def\O{{\sy T}}
\def\C#1 #2{{\(#1\)_{#2}}}

\def\D{\L}

\def\ttchap#1{
\numlem=0
\vfill\eject\gauche={\the\numchap.
#1}\droite={\hfill}{\bfXVIII \BX{CHAPITRE 
\the\numchap. \bn \hbox{#1}} }\bn\vskip 2cm}

\def\ttparag
#1{{\bfXII #1}\bn}

 \def\Pr{\bg\bf PREUVE.
\rm} \def\Rp{\hbox{\carre\rm}\bn}

\def\Lem#1#2{\advance\numlem by
1\bg{{\bf LEMME \the\numlem. #1}\sl#2}\rm\bg}
\def\Lemprop#1#2{\advance\numlem by 1\bg{{\bf
PROPOSITION \the\numlem. #1}\sl#2}\rm\bg}
\def\Lempropr#1#2{\advance\numlem by 1\bg{{\bf
PROPRIETE \the\numlem. #1}\sl#2}\rm\bg}
\def\Thm#1#2{\advance\numlem by 1\bg{{\bf
THEOREME \the\numlem. #1}\sl#2}\rm\bg}

\font\sc=cmcsc10

\def\ttchap#1{
\numlem=0
\vfill\eject\gauche={\the\numchap.
#1}\droite={\hfill}{\bfXVIII \BX{CHAPTER  \the\numchap.\bn
\hbox{#1}} }\bn\vskip 2cm}

\def\Pr{\bg\bf Proof.
\rm}

\def\Lem#1#2{\advance\numlem by
1\bg{{\bf Lemma \the\numlem. #1}\sl#2}\rm\bg}
\def\Lemprop#1#2{\advance\numlem by 1\bg{{\bf
Proposition \the\numlem. #1}\sl#2}\rm\bg}
\def\Prop#1#2{\bg{{\bf
Proposition  #1}\sl#2}\rm\bg}
\def\Lempropr#1#2{\advance\numlem by 1\bg{{\bf
Property \the\numlem. #1}\sl#2}\rm\bg}
\def\Thm#1#2{\advance\numlem by 1\bg{{\bf
Theorem \the\numlem. #1}\sl#2}\rm\bg}

\catcode`\Ž=\active \def Ž{\'e}
\catcode`\ˆ=\active \def ˆ{\`a}
\catcode`\=\active \def {\`u}
\catcode`\=\active \def {\`e}
\catcode`\=\active \def {\^e}
\catcode`\™=\active \def ™{\^o}
\catcode`\"=\active \def "{\^i}
\catcode`\=\active \def {\c c}

\def\fileversion{v1.0}
\def\filedate{3 Dec 91}
\immediate\write16{Document style option `epsfig', \fileversion\space
<\filedate> (edited by SPQR)}
\chardef\atcode=\catcode`\@
\catcode`\@=11\relax
\ifx\typeout\undefined%
        \newwrite\@unused
        \def\typeout#1{{\let\protect\string\immediate\write\@unused{#1}}}
\fi
\ifx\undefined\epsfig%
\else
        \typeout{EPSFIG --- already loaded} 
\fi
%
%
\ifx\undefined\fbox\def\fbox#1{#1}\fi
%
\ifx\undefined\epsfbox\input epsf\fi
%
\ifx\undefined\@latexerr
        \newlinechar`\^^J
        \def\@spaces{\space\space\space\space}
        \def\@latexerr#1#2{%
        \edef\@tempc{#2}\expandafter\errhelp\expandafter{\@tempc}%
        \typeout{Error. \space see a manual for explanation.^^J
         \space\@spaces\@spaces\@spaces Type \space H <return> \space for
         immediate help.}\errmessage{#1}}
\fi
\def\@whattodo{You tried to include a PostScript figure which 
cannot be found^^JIf you press return to carry on anyway,^^J
The failed name will be printed in place of the figure.^^J
or type X to quit}
\def\@whattodobb{You tried to include a PostScript figure which 
has no^^Jbounding box, and you supplied none.^^J
If you press return to carry on anyway,^^J
The failed name will be printed in place of the figure.^^J
or type X to quit}
%
\newwrite\@unused
%
\def\@nnil{\@nil}
\def\@empty{}
\def\@psdonoop#1\@@#2#3{}
\def\@psdo#1:=#2\do#3{\edef\@psdotmp{#2}\ifx\@psdotmp\@empty \else
    \expandafter\@psdoloop#2,\@nil,\@nil\@@#1{#3}\fi}
\def\@psdoloop#1,#2,#3\@@#4#5{\def#4{#1}\ifx #4\@nnil \else
       #5\def#4{#2}\ifx #4\@nnil \else#5\@ipsdoloop #3\@@#4{#5}\fi\fi}
\def\@ipsdoloop#1,#2\@@#3#4{\def#3{#1}\ifx #3\@nnil 
       \let\@nextwhile=\@psdonoop \else
      #4\relax\let\@nextwhile=\@ipsdoloop\fi\@nextwhile#2\@@#3{#4}}
\def\@tpsdo#1:=#2\do#3{\xdef\@psdotmp{#2}\ifx\@psdotmp\@empty \else
    \@tpsdoloop#2\@nil\@nil\@@#1{#3}\fi}
\def\@tpsdoloop#1#2\@@#3#4{\def#3{#1}\ifx #3\@nnil 
       \let\@nextwhile=\@psdonoop \else
      #4\relax\let\@nextwhile=\@tpsdoloop\fi\@nextwhile#2\@@#3{#4}}
%
%
%
\long\def\epsfaux#1#2:#3\\{\ifx#1\epsfpercent
   \def\testit{#2}\ifx\testit\epsfbblit
        \@atendfalse
        \epsf@atend #3 . \\%
        \if@atend       
           \if@verbose
                \typeout{epsfig: found `(atend)'; continuing search}
           \fi
        \else
                \epsfgrab #3 . . . \\%
                \global\no@bbfalse
        \fi
   \fi\fi}%
%
%
\def\epsf@atendlit{(atend)} 
\def\epsf@atend #1 #2 #3\\{%
   \def\epsf@tmp{#1}\ifx\epsf@tmp\empty
      \epsf@atend #2 #3 .\\\else
   \ifx\epsf@tmp\epsf@atendlit\@atendtrue\fi\fi}


\chardef\letter = 11
\chardef\other = 12

\newif \ifdebug 
\newif\ifc@mpute 
\newif\if@atend
\c@mputetrue 

\let\then = \relax
\def\r@dian{pt }
\let\r@dians = \r@dian
\let\dimensionless@nit = \r@dian
\let\dimensionless@nits = \dimensionless@nit
\def\internal@nit{sp }
\let\internal@nits = \internal@nit
\newif\ifstillc@nverging
\def \Mess@ge #1{\ifdebug \then \message {#1} \fi}

{ 
        \catcode `\@ = \letter
        \gdef \nodimen {\expandafter \n@dimen \the \dimen}
        \gdef \term #1 #2 #3%
               {\edef \t@ {\the #1}
                \edef \t@@ {\expandafter \n@dimen \the #2\r@dian}%
                \t@rm {\t@} {\t@@} {#3}%
               }
        \gdef \t@rm #1 #2 #3%
               {{%
                \count 0 = 0
                \dimen 0 = 1 \dimensionless@nit
                \dimen 2 = #2\relax
                \Mess@ge {Calculating term #1 of \nodimen 2}%
                \loop
                \ifnum  \count 0 < #1
                \then   \advance \count 0 by 1
                        \Mess@ge {Iteration \the \count 0 \space}%
                        \Multiply \dimen 0 by {\dimen 2}%
                        \Mess@ge {After multiplication, term = \nodimen 0}%
                        \Divide \dimen 0 by {\count 0}%
                        \Mess@ge {After division, term = \nodimen 0}%
                \repeat
                \Mess@ge {Final value for term #1 of 
                                \nodimen 2 \space is \nodimen 0}%
                \xdef \Term {#3 = \nodimen 0 \r@dians}%
                \aftergroup \Term
               }}
        \catcode `\p = \other
        \catcode `\t = \other
        \gdef \n@dimen #1pt{#1} 
}

\def \Divide #1by #2{\divide #1 by #2} 

\def \Multiply #1by #2
       {{
        \count 0 = #1\relax
        \count 2 = #2\relax
        \count 4 = 65536
        \Mess@ge {Before scaling, count 0 = \the \count 0 \space and
                        count 2 = \the \count 2}%
        \ifnum  \count 0 > 32767 
        \then   \divide \count 0 by 4
                \divide \count 4 by 4
        \else   \ifnum  \count 0 < -32767
                \then   \divide \count 0 by 4
                        \divide \count 4 by 4
                \else
                \fi
        \fi
        \ifnum  \count 2 > 32767 
        \then   \divide \count 2 by 4
                \divide \count 4 by 4
        \else   \ifnum  \count 2 < -32767
                \then   \divide \count 2 by 4
                        \divide \count 4 by 4
                \else
                \fi
        \fi
        \multiply \count 0 by \count 2
        \divide \count 0 by \count 4
        \xdef \product {#1 = \the \count 0 \internal@nits}%
        \aftergroup \product
       }}

\def\r@duce{\ifdim\dimen0 > 90\r@dian \then   
                \multiply\dimen0 by -1
                \advance\dimen0 by 180\r@dian
                \r@duce
            \else \ifdim\dimen0 < -90\r@dian \then  
                \advance\dimen0 by 360\r@dian
                \r@duce
                \fi
            \fi}

\def\Sine#1%
       {{%
        \dimen 0 = #1 \r@dian
        \r@duce
        \ifdim\dimen0 = -90\r@dian \then
           \dimen4 = -1\r@dian
           \c@mputefalse
        \fi
        \ifdim\dimen0 = 90\r@dian \then
           \dimen4 = 1\r@dian
           \c@mputefalse
        \fi
        \ifdim\dimen0 = 0\r@dian \then
           \dimen4 = 0\r@dian
           \c@mputefalse
        \fi
        \ifc@mpute \then
                \divide\dimen0 by 180
                \dimen0=3.141592654\dimen0
                \dimen 2 = 3.1415926535897963\r@dian 
                \divide\dimen 2 by 2 
                \Mess@ge {Sin: calculating Sin of \nodimen 0}%
                \count 0 = 1 
                \dimen 2 = 1 \r@dian 
                \dimen 4 = 0 \r@dian 
                \loop
                        \ifnum  \dimen 2 = 0 
                        \then   \stillc@nvergingfalse 
                        \else   \stillc@nvergingtrue
                        \fi
                        \ifstillc@nverging 
                        \then   \term {\count 0} {\dimen 0} {\dimen 2}%
                                \advance \count 0 by 2
                                \count 2 = \count 0
                                \divide \count 2 by 2
                                \ifodd  \count 2 
                                \then   \advance \dimen 4 by \dimen 2
                                \else   \advance \dimen 4 by -\dimen 2
                                \fi
                \repeat
        \fi             
                        \xdef \sine {\nodimen 4}%
       }}

\def\Cosine#1{\ifx\sine\UnDefined\edef\Savesine{\relax}\else
                             \edef\Savesine{\sine}\fi
        {\dimen0=#1\r@dian\multiply\dimen0 by -1
         \advance\dimen0 by 90\r@dian
         \Sine{\nodimen 0}
         \xdef\cosine{\sine}
         \xdef\sine{\Savesine}}}              
%
\def\psdraft{\def\@psdraft{0}}
\def\psfull{\def\@psdraft{100}}
\psfull
\newif\if@draftbox
\def\psnodraftbox{\@draftboxfalse}
\@draftboxtrue
\newif\if@noisy
\def\pssilent{\@noisyfalse}
\def\psnoisy{\@noisytrue}
\@noisyfalse

\newif\if@bbllx
\newif\if@bblly
\newif\if@bburx
\newif\if@bbury
\newif\if@height
\newif\if@width
\newif\if@rheight
\newif\if@rwidth
\newif\if@angle
\newif\if@clip
\newif\if@verbose
\def\@p@@sclip#1{\@cliptrue}
\newif\if@prologfile
\@prologfiletrue
%
\newif\ifuse@psfig
\def\@p@@sfile#1{%
\def\@p@sfile{NO FILE: #1}%
\def\@p@sfilefinal{NO FILE: #1}%
        \openin1=#1
        \ifeof1\closein1
                \openin1=#1.bb
                        \ifeof1\closein1
                                \if@bbllx\if@bblly\if@bburx\if@bbury
                                        \def\@p@sfile{#1}%
                                        \def\@p@sfilefinal{#1}%
                                        \fi\fi\fi
                                \else
                                        \@latexerr{ERROR! PostScript file #1 not found}\@whattodo
                                        \@p@@sbbllx{100bp}
                                        \@p@@sbblly{100bp}
                                        \@p@@sbburx{200bp}
                                        \@p@@sbbury{200bp}
                                        \def\@p@scost{200}
                                \fi
                        \else
                                \closein1%
                                \edef\@p@sfile{#1.bb}%
                                \edef\@p@sfilefinal{"`zcat `texfind #1.Z`"}%
                        \fi
        \else\closein1
                    \edef\@p@sfile{#1}%
                    \edef\@p@sfilefinal{#1}%
        \fi%
}
\let\@p@@sfigure\@p@@sfile
\def\@p@@sbbllx#1{
                \use@psfigtrue
                \@bbllxtrue
                \dimen100=#1
                \edef\@p@sbbllx{\number\dimen100}
}
\def\@p@@sbblly#1{
                \use@psfigtrue
                \@bbllytrue
                \dimen100=#1
                \edef\@p@sbblly{\number\dimen100}
}
\def\@p@@sbburx#1{
                \use@psfigtrue
                \@bburxtrue
                \dimen100=#1
                \edef\@p@sbburx{\number\dimen100}
}
\def\@p@@sbbury#1{
                \use@psfigtrue
                \@bburytrue
                \dimen100=#1
                \edef\@p@sbbury{\number\dimen100}
}
\def\@p@@sheight#1{
                \@heighttrue
                \epsfysize=#1
                \dimen100=#1
                \edef\@p@sheight{\number\dimen100}
}
\def\@p@@swidth#1{
                \@widthtrue
                \epsfxsize=#1
                \dimen100=#1
                \edef\@p@swidth{\number\dimen100}
}
\def\@p@@srheight#1{
                \@rheighttrue
                \dimen100=#1
                \edef\@p@srheight{\number\dimen100}
}
\def\@p@@srwidth#1{
                \@rwidthtrue
                \dimen100=#1
                \edef\@p@srwidth{\number\dimen100}
}
\def\@p@@sangle#1{
                \use@psfigtrue
                \@angletrue
                \edef\@p@sangle{#1} 
}
\def\@p@@ssilent#1{ 
                \@verbosefalse
}
\def\@cs@name#1{\csname #1\endcsname}
\def\@setparms#1=#2,{\@cs@name{@p@@s#1}{#2}}
%
%
\def\ps@init@parms{
                \@bbllxfalse \@bbllyfalse
                \@bburxfalse \@bburyfalse
                \@heightfalse \@widthfalse
                \@rheightfalse \@rwidthfalse
                \def\@p@sbbllx{}\def\@p@sbblly{}
                \def\@p@sbburx{}\def\@p@sbbury{}
                \def\@p@sheight{}\def\@p@swidth{}
                \def\@p@srheight{}\def\@p@srwidth{}
                \def\@p@sangle{0}
                \def\@p@sfile{}
                \def\@p@scost{10}
                \use@psfigfalse
                \def\@sc{}
                \if@noisy
                        \@verbosetrue
                \else
                        \@verbosefalse
                \fi
                \@clipfalse
}
%
%
\def\parse@ps@parms#1{
                \@psdo\@psfiga:=#1\do
                   {\expandafter\@setparms\@psfiga,}}
%
%
\newif\ifno@bb
\def\bb@missing{
        \epsfgetbb{\@p@sfile}
        \ifepsfbbfound\no@bbfalse\else\no@bbtrue\bb@cull\epsfllx\epsflly\epsfurx\epsfury\fi
}       
\def\bb@cull#1#2#3#4{
        \dimen100=#1 bp\edef\@p@sbbllx{\number\dimen100}
        \dimen100=#2 bp\edef\@p@sbblly{\number\dimen100}
        \dimen100=#3 bp\edef\@p@sbburx{\number\dimen100}
        \dimen100=#4 bp\edef\@p@sbbury{\number\dimen100}
        \no@bbfalse
}

\newdimen\p@intvaluex
\newdimen\p@intvaluey
\newdimen\@ffsetvalue
\newdimen\x@ffsetvalue
\newdimen\y@ffsetvalue


\def\compute@offset#1#2{{\dimen0=#1 sp\dimen1=#2 sp
                        \advance\dimen1 by -\dimen0
                        \dimen1=\sine\dimen1
                        \dimen0=\cosine\dimen1
                        \ifdim\dimen0<0sp \dimen1=0sp \fi
                        \global\@ffsetvalue=\dimen1}}

\def\rotate@#1#2{{\dimen0=#1 sp\dimen1=#2 sp
                  \global\p@intvaluex=\cosine\dimen0
                  \dimen3=\sine\dimen1
                  \global\advance\p@intvaluex by -\dimen3
                  \global\p@intvaluey=\sine\dimen0
                  \dimen3=\cosine\dimen1
                  \global\advance\p@intvaluey by \dimen3
                  }}
%
\def\compute@bb{
                \no@bbfalse
                \if@bbllx \else \no@bbtrue \fi
                \if@bblly \else \no@bbtrue \fi
                \if@bburx \else \no@bbtrue \fi
                \if@bbury \else \no@bbtrue \fi
                \ifno@bb \bb@missing \fi
                \ifno@bb
                        \@latexerr{ERROR! cannot locate BB!}\@whattodobb
                        \@p@@sbbllx{100bp}
                        \@p@@sbblly{100bp}
                        \@p@@sbburx{200bp}
                        \@p@@sbbury{200bp}
                        \def\@p@scost{200}
                \fi
                \if@angle 
                        \Sine{\@p@sangle}\Cosine{\@p@sangle}
                        \compute@offset{\@p@sbblly}{\@p@sbbury}
                        \x@ffsetvalue=\@ffsetvalue
                        \compute@offset{\@p@sbburx}{\@p@sbbllx}
                        \y@ffsetvalue=\@ffsetvalue

                        \rotate@{\@p@sbbllx}{\@p@sbblly}
                        \advance\p@intvaluex by -\x@ffsetvalue
                        \advance\p@intvaluey by -\y@ffsetvalue
                        \edef\@p@sbbllx{\number\p@intvaluex}
                        \edef\@p@sbblly{\number\p@intvaluey}

                        \rotate@{\@p@sbburx}{\@p@sbbury}
                        \advance\p@intvaluex by \x@ffsetvalue
                        \advance\p@intvaluey by \y@ffsetvalue
                        \edef\@p@sbburx{\number\p@intvaluex}
                        \edef\@p@sbbury{\number\p@intvaluey}
                        {
                         \count0=\@p@sbbllx \count1=\@p@sbblly
                         \count2=\@p@sbburx \count3=\@p@sbbury
                         \dimen0=\@p@sbbllx sp\dimen1=\@p@sbblly sp
                         \dimen2=\@p@sbburx sp\dimen3=\@p@sbbury sp
                         \dimen203=\dimen2 \advance\dimen203 by -\dimen0
                         \dimen204=\dimen3 \advance\dimen204 by -\dimen1
                         \ifdim\dimen203<0sp 
                              \count203=\count2 \count2=\count0 
                              \count0=\count203 
                              \global\edef\@p@sbbllx{\number\count0}
                              \global\edef\@p@sbburx{\number\count2}
                         \fi
                         \ifdim\dimen204<0sp 
                               \count204=\count3
                               \count3=\count1
                               \count1=\count204
                               \global\edef\@p@sbblly{\number\count1}
                               \global\edef\@p@sbbury{\number\count3}
                         \fi
                        }
                \fi
                \count203=\@p@sbburx
                \count204=\@p@sbbury
                \advance\count203 by -\@p@sbbllx
                \advance\count204 by -\@p@sbblly
                \edef\@bbw{\number\count203}
                \edef\@bbh{\number\count204}
}
%
%
\def\in@hundreds#1#2#3{\count240=#2 \count241=#3
                     \count100=\count240        
                     \divide\count100 by \count241
                     \count101=\count100
                     \multiply\count101 by \count241
                     \advance\count240 by -\count101
                     \multiply\count240 by 10
                     \count101=\count240        
                     \divide\count101 by \count241
                     \count102=\count101
                     \multiply\count102 by \count241
                     \advance\count240 by -\count102
                     \multiply\count240 by 10
                     \count102=\count240        
                     \divide\count102 by \count241
                     \count200=#1\count205=0
                     \count201=\count200
                        \multiply\count201 by \count100
                        \advance\count205 by \count201
                     \count201=\count200
                        \divide\count201 by 10
                        \multiply\count201 by \count101
                        \advance\count205 by \count201
                     \count201=\count200
                        \divide\count201 by 100
                        \multiply\count201 by \count102
                        \advance\count205 by \count201
                     \edef\@result{\number\count205}
}
\def\compute@wfromh{
                \in@hundreds{\@p@sheight}{\@bbw}{\@bbh}
                \edef\@p@swidth{\@result}
}
\def\compute@hfromw{
                \in@hundreds{\@p@swidth}{\@bbh}{\@bbw}
                \edef\@p@sheight{\@result}
}
\def\compute@handw{
                \if@height 
                        \if@width
                        \else
                                \compute@wfromh
                        \fi
                \else 
                        \if@width
                                \compute@hfromw
                        \else
                                \edef\@p@sheight{\@bbh}
                                \edef\@p@swidth{\@bbw}
                        \fi
                \fi
}
\def\compute@resv{
                \if@rheight \else \edef\@p@srheight{\@p@sheight} \fi
                \if@rwidth \else \edef\@p@srwidth{\@p@swidth} \fi
}
%
\def\compute@sizes{
        \compute@bb
        \compute@handw
        \compute@resv
}
%
%
\def\psfig#1{\vbox {
        %
        \ps@init@parms
        \parse@ps@parms{#1}
        \ifnum\@p@scost<\@psdraft
                \typeout{[\@p@sfilefinal]}
                \if@verbose
                        \typeout{epsfig: using PSFIG macros}
                \fi
                \psfig@method
        \else
                \epsfig@draft
        \fi
}}

\def\epsfig#1{\vbox {
        %
        \ps@init@parms
        \parse@ps@parms{#1}
        \ifnum\@p@scost<\@psdraft
                \typeout{[\@p@sfilefinal]}
                \if@clip\use@psfigtrue\fi
                \if@angle\use@psfigtrue\fi
                \ifuse@psfig
                        \if@verbose
                                \typeout{epsfig: using PSFIG macros}
                        \fi
                        \psfig@method
                \else
                        \if@verbose
                                \typeout{epsfig: using EPSF macros}
                        \fi
                        \epsf@method
                \fi
        \else
                \epsfig@draft
        \fi
}}

\def\epsf@method{%
        \epsfgetbb{\@p@sfile}%
        \epsfsetgraph{\@p@sfilefinal}
}
\def\psfig@method{%
        \compute@sizes
        \special{ps::[begin]  \@p@swidth \space \@p@sheight \space%
        \@p@sbbllx \space \@p@sbblly \space%
        \@p@sbburx \space \@p@sbbury \space%
        startTexFig \space }%
        \if@angle
                \special {ps:: \@p@sangle \space rotate \space} 
        \fi
        \if@clip
                \if@verbose
                        \typeout{(clipped to BB) }
                \fi
                \special{ps:: doclip \space }%
        \fi
        \special{ps: plotfile \@p@sfilefinal \space }%
        \special{ps::[end] endTexFig \space }%
        \vbox to \@p@srheight true sp{\hbox to \@p@srwidth true sp{\hss}\vss}
}
%
\def\epsfig@draft{
\compute@sizes
\if@draftbox
        \hbox{\fbox{\vbox to \@p@srheight true sp{
        \vss\hbox to \@p@srwidth true sp{ \hss \@p@sfilefinal \hss }\vss
        }}}
\else
        \vbox to \@p@srheight true sp{%
        \vss\hbox to \@p@srwidth true sp{\hss}\vss}
\fi     
}
\catcode`\@=\atcode 

\newcount\dcoef
\newcount\mcoef

\def\dessin(#1)(#2)(#3){{\epsfig{file=#1,width=#2truecm,height=#3truecm}}}

\mcoef=1
\dcoef=2

\gauche={S. Burckel}
\droite={Complexity of some Path Problems in DAGs and Linear Orders.}

\centerline{\XIV Complexity of some Path Problems in}
\centerline{\XIV DAGs and Linear Orders.}
\sk

\bn {\bf Abstract. We investigate here the computational complexity of three natural problems in directed acyclic graphs.
We prove their NP Completeness and consider their restrictions to linear orders.
}

\bn {\it Mathematics Subjects Classification : .}

\bn
\ttparag{1. Introduction.}

In this paper, we precise and extend the results of [1] about the complexity of path problems in directed graphs.
In particular we focus ourselves on the case of DAGs and very particular ones : linear orders.

\bn
\ttparag{2. Problems.}

\bn {\bf Problem 1 : Null Weighted Path.}
\bn
input : a DAG $G=(V,E)$ endowed with a weight function $w:E\mapsto \ZZ$ and two vertices $s,t$.

question : is there a directed path $\pi$ from $s$ to $t$ such that 
$$\sum_{e\in\pi}w(e)=0$$

\Prop{(1.)}{ The Null Weighted Path Decision Problem is NP-Complete, even on linear orders.}

\Pr
Obviously, it is in NP. To prove that it is NP hard, we reduce the SUBSET SUM Problem to it. Let $A=\{a_1,\dots,a_n\}$ be a subset of $\ZZ^n$.
One can construct in polynomial time a DAG $G=(V,E)$, a weight function and two vertices $s,t$ that admits a path from $s$ to $t$ of null weight
if and only if there exists a non empty set $S\subseteq A$ such that $$\sum_{a\in S}a=0$$
 
The DAG $G=(V,E)$ has $n+2$ ordered vertices $V=[s=0,1,2,\dots,n,n+1=t]$ and the arcs are all the pairs $e_{ij}=(i,j)$ with $i<j$. Hence, $G$ is a linear order.
The weights are~:
$$w(e_{ij})=\cases{
a_j& for $0\le i<j\le n$\cr
0& for $1\le i\le n$ and $j=t$\cr
+1& for $i=0$ and $j=t$
}$$
\sk The following example shows this construction for $A=\{4,2,-5\}$ (arcs are oriented from left to right)~:

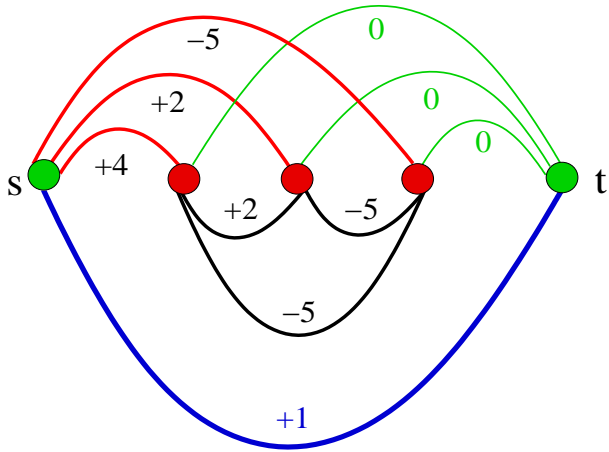
$$\dessin(pb1.eps)(8)(6)$$

\sk Assume there exists a non empty subset of $A$ $S=\{a_{i_1},\dots,a_{i_k}\}$ such that $a_{i_1}+\dots+a_{i_k}=0$. Then,
one can assume that $i_1<i_2<\dots<i_k$. Hence, 
$\pi=[s,a_{i_1},\dots,a_{i_k},t]$ is a null weighted path from $s$ to $t$.
\sk For the converse, observe that every path from $s$ to $t$ (excepted the direct arc $(s,t)$ of weight $+1$)  
contains at least one arc with an element of $A$. Hence, a null weighted path from $s$ to $t$ describes a non empty subset of $A$ with null sum.
\Rp

\bn
Now we consider a similar problem on DAGs with positive weights.

\bn {\bf Problem 2 : $K$ Weighted Path.}
\bn
input : a DAG $G=(V,E)$ endowed with a weight function $w:E\mapsto \NN$, two vertices $s,t$ and an integer $K\ge 0$

question : is there a directed path $\pi$ from $s$ to $t$ such that 
$$\sum_{e\in\pi}w(e)=K$$

\Prop{(2.)}{ The $K$ Weighted Path Decision Problem is NP-Complete, even on linear orders.}

\Pr
Obviously, it is in NP. To prove that it is NP hard, we reduce again the SUBSET SUM Problem to it. Let $A=\{a_1,\dots,a_n\}$ be a subset of $\ZZ^n$.
One can construct in polynomial time a DAG $G=(V,E)$, a weight function, an integer $K$ and two vertices $s,t$ that admits a path from $s$ to $t$ of weight $K$ 
if and only if there exists a non empty set $S\subseteq A$ such that $$\sum_{a\in S}a=0$$ 
\sk Let $P=-min(\{0\}\cup A)$ and fix $K=(n+1)P$.
\sk
The construction is similar as previously. However, one translates the weight of each arc $(i,j)$ by $(j-i)P$.
The DAG $G=(V,E)$ has $n+2$ ordered vertices $V=[s=0,1,2,\dots,n,n+1=t]$ and the arcs are all the pairs $e_{ij}=(i,j)$ with $i<j$. Hence $G$ is a linear order.
The weights are~:
$$w(e_{ij})=\cases{
(j-i)P+a_j& for $0\le i<j\le n$\cr
(j-i)P& for $1\le i\le n$ and $j=t$\cr
(j-i)P+1& for $i=0$ and $j=t$
}$$
\sk For $A=\{4,2,-5\}$ one obtains with $P=-min(\{0,4,2,-5\})=+5$ the following linear order DAG~:

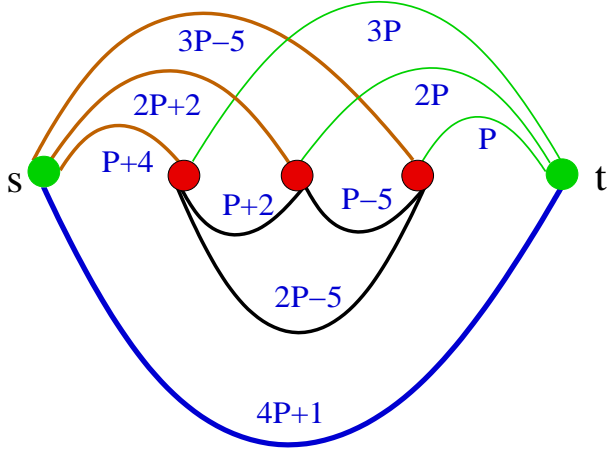
$$\dessin(pb2.eps)(8)(6)$$

\sk Assume there exists a non empty subset of $A$ $S=\{a_{i_1},\dots,a_{i_k}\}$ such that $a_{i_1}+\dots+a_{i_k}=0$. Then,
one can assume that $i_1<i_2<\dots<i_k$. Hence, 
$\pi=[s,a_{i_1},\dots,a_{i_k},t]$ is a path from $s$ to $t$ of weight $K=(n+1)P$.
\sk For the converse, observe that every path from $s$ to $t$ (excepted the direct arc $(s,t)$ of weight $(n+1)P+1$) contains at least one arc that contributes to 
an element of $A$. Hence, a path of weight $(n+1)P$ describes a non empty subset of $A$ with null sum.
\Rp

\bn
Now, we consider DAGs without weights. We just count the length of the paths, i.e, the number of arcs in them.

\bn {\bf Problem : Path of Length $K$.}
\bn
input : a DAG $G=(V,E)$, two vertices $s,t$ and an integer $K>0$.

question : is there a directed path from $s$ to $t$ of length $K$.

\bn
This problem is solvable in deterministic polynomial time since one can assume that $K<|V|$ since $G$ is a DAG. Hence, one can 
efficiently compute the matrix $M^K$ where $M$ is the adjacency matrix of $G$.
\bn
However, when one adds a supplementary condition, the problem becomes NP-Complete.

\bn {\bf Problem 3 : Irreducible Path of Length $K$.}
\bn
input : a DAG $G=(V,E)$, two vertices $s,t$ and an integer $K>0$.

question : is there an irreducible directed path from $s$ to $t$ of length $K$.
\bn
Here {\it irreducible path} means a sequence of vertices $[s=x_0,x_1,\dots,x_K=t]$ such that $(x_i,x_j)\in E$ if and only if $j=i+1$.
Of course, in linear orders, this problem is trivial since~:
\sk for $K\ge 2$, the answer is necessarily NO (by transitivity).
\sk for $K=1$, that corresponds to $s<t$. 
\bn
For more general DAGs, this problem is non trivial anymore.

\Prop{(3.)}{ The Irreducible Path of Length $K$ Decision Problem is NP-Complete.}

\Pr
Obviously again, it is in NP. For the NP Hardness, we reduce now the CNF SAT Problem to it.
Let $\Phi$ be a CNF formula $C_1\wedge C_2\wedge C_3\dots \wedge C_k$ where each clause $C_i$ is a disjunction of litterals
on different variables. Assume that the total number of litterals in $\Phi$ is $N$. 
We construct a DAG $G=(V,E)$ with $N+2$ vertices~: one vertex per litteral in each clause $C_i$ and two supplementary vertices $s,t$.
The arcs are the pairs~:
\sk $(s,x)$ for $x\in C_1$
\sk $(y,t)$ for $y\in C_k$
\sk $(x,y)$ for $x\in C_i$ and $y\in C_{j}$ with $j=i+1$ and $x\not=\bar{y}$
\sk $(x,y)$ for $x\in C_i$ and $y\in C_{j}$ with $j>i+1$ and $x=\bar{y}$
\bn
For $\Phi=(a\vee b\vee \bar{c})\wedge(d\vee \bar{a})\wedge(\bar{a}\vee \bar{d}\vee c)$, one obtains the following DAG~:

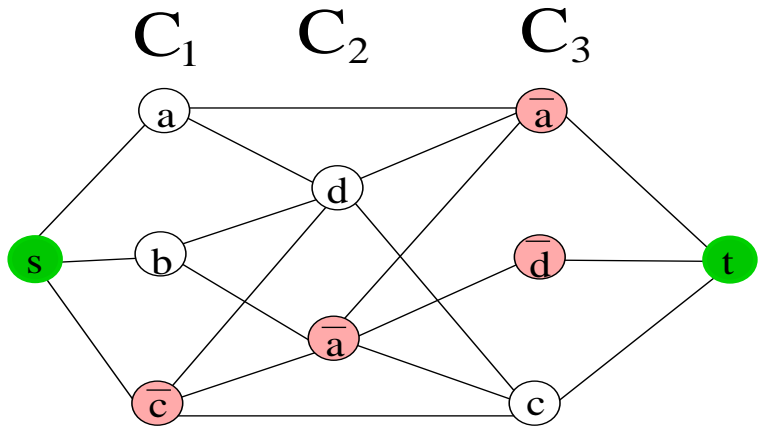
$$\dessin(pb3.eps)(10)(6)$$

The question is to find an irreducible path of length $K=(k+1)$.
This is equivalent to the satisfiability of $\Phi$~: one has to find at least one litteral assigned to "true" in each clause (hence a path of length $k+1$ from
$s$ to $t$) and with no contradictory assignments like $x=true$ and $\bar{x}=true$. 
For consecutive clauses, such an arc $(x,\bar{x})$ does not exist. Otherwise, this arc $(x,\bar{x})$
would make the path not irreducible.

\Rp

\bn


\def \it#1{#1}

\parindent=1.8truecm \def\Ref#1; #2; #3; #4;{ \sgni\item{\hbox to 1.6truecm{ [#1]\hfill}}{\sc #2},
{\sl #3}, #4. } \def\Reff#1; #2; #3; #4; #5; #6; #7;{ \sgni\item{\hbox to 1.6truecm{
[#1]\hfill}}{\sc #2}, {\sl #3}, #4 {\bf #5} (#6) #7. }

\Reff {1}; S. Basagni,D. Bruschi, F. Ravasio; On the difficulty of finding walks of length $k$; Theoretical Informatics and Applications; 31 (5); 1997;
429--435;

\bn

\sk\ni Serge Burckel. 
\sk\ni INRIA-LORIA,
\sk\ni Campus Scientifique
\sk\ni BP 239
\sk\ni 54506 Vandoeuvre-l\`es-Nancy Cedex 
\sk\ni sergeburckel@orange.fr

\bn
{\it This work has been partially supported by the CNRS-ANR project "graphs decompositions and algorithms".}

\end